Discussion

# Unity and Disunity in Mathematics

Bernhelm Booß-Bavnbek (Roskilde University, Denmark) and Philip J. Davis (Providence, RI, USA)

**Foreword**
While dictionaries have produced a number of definitions of the word "unity", the sense in which we employ the terms unity or disunity will emerge from what follows. The unity of mathematics has been proclaimed and celebrated by numerous mathematicians of the first rank, most recently by M. F. Atiyah and I. M. Gelfand. Atiyah has found the notion of unity in the interaction between geometry and physics while Gelfand asserts that its "beauty, simplicity [!], exactness, and crazy ideas" as well as its capability of abstraction and generalisation are all hallmarks of an inner unity.

In other places, we have elaborated this by pointing out that the unity of mathematics has, as a consequence, its power to compactify experiences in a form capable of being transferred and modified or adapted to new mathematical situations. Existing knowledge is mustered towards creating the iterative sequence description, prescription/design, execution, comparison.

An additional claim of unity is that mathematics research for at least the last 200 years has been largely dominated by a single methodology: that of rigorous proof. These assertions add up to the "Unity of Mathematics", a well-known and well publicised, predominantly philosophical claim that has been promoted in international conferences on the topic.

Important declarations can be found in the volume *The Unity of Mathematics* edited by Etingof et al. [7]. In this book about 20 articles on a variety of topics, written by preeminent authorities, have been anthologised and advanced under the slogan of unity. What a casual reader may wonder about this collection is whether the various authors are capable of or are interested in understanding what their fellow authors have written. The second named author, a professional who has functioned for some years as a researcher, an applied mathematician and an educator, found all the articles incomprehensible. He rejects the argument that since "mathematics is a universal language", with a bit of study he would find them all understandable. This would be similar to the assertion that with a bit of goodwill in the universe, the lion will lie down with the lamb and universal peace will reign on earth. Thus, while the editors trumpet the unity of mathematics, the contents of the book assert equally well its disunity. This is an instance of a simultaneous and dual nature, expressed by unity/disunity, that permeates the field.

Grattan-Guinness ([8], p. 490) points out that by the turn of the 20th century, "a snobbish preference for pure mathematics over applications became more marked … so that the profession tended to separate into two components. In Berlin and elsewhere in Germany, purism became an explicit creed for the professional…"

We have had (and still have) the *Journal für die Reine und Angewandte Mathematik* (Crelle's Journal, founded 1826)[1] and the *Journal de Mathématiques Pures et Appliquées* (Liouville's Journal, founded 1836). Thus, while the existence of what is pure and what is applied would seem to proclaim disunity, they often stand united.

Yet, we believe that the phrase "unity of mathematics" expresses a dream, an ideal that doesn't exist. We shall point to diachronic and cross cultural disunities, to semantic, semiotic and philosophic ambiguities and to the non-acceptance of certain mathematical texts by some practitioners of the subject. But to deny the strivings towards unity would be misleading and not constructive.

**Extent and practices of the mathematical enterprise**
The richness of the mathematical enterprise can be exhibited in many different ways and can be used to show both its inherent unity and its disunity. The number of mathematicians currently in the world has been estimated at 100,000. The number of mathematics books in the Brown University Science Library is approximately 32,000. Many of these books carry a title that implies the gathering together and the wedding of information around a single topic. Within this "marriage" unity may be found, e.g. *Analytic Function Theory*, *Circulant Matrices* and *Elliptic Boundary Problems for Dirac Operators*.

Such books display the mathematical enterprise at a high level of theory and application. At the informal level of discourse and practice, a life today without the continuous use of mathematics in commerce, finance, science and technology, medical practice, art, sports, etc. is inconceivable. We are living in a highly mathematised social and technological world.

Let us look at the published material. The subject classification scheme of the Mathematical Reviews and Zentralblatt MATH (MSC 2010) at the macro level lists almost a *hundred* mathematical subjects, random examples being: 00 General; 13 Commutative Rings and Algebras; 42 Fourier Analysis; 83 Relativity and Gravitational Theory; and 97 Mathematical Education. At the micro level, within each group there are subgroups such as 13A15: Ideals and Multiplicative Ideal Theory.

Since mathematics has the aspect of a language, it is instructive to compare it with the classification schemes for literature and rhetoric. The Dewey Decimal Classification System (1878 and with constant updates) at the macro level lists 100 categories numbered from 800 to 899. Since mathematics also has the aspect of an art, it would be of interest to compare its classification schemes with the MSC 2010.

---

[1] In an old and sceptical play on words, Reine und Angewandte Mathematik become Reine und*un*angewandte Mathematik (pure inapplicable mathematics). Notice also the predominant plural in the French use of the term *mathématique*.

To some extent, each macro and micro mathematical subject has its own techniques and its own intellectual resources and devotees who both utilise and enrich its resources. While there may, indeed, be some connections between e.g. potential theory and non-associative algebras, elaborate abstract graphs and tree structures have been drawn, linking or uniting selected mathematical topics. This can be done even with an individual theorem. The second named author recalls that as a graduate student in Professor Hassler Whitney's course on topology, he made a tree structure for one of the major theorems in topological group theory showing how many significant theorems fed into the final statement as sub-results.

The frequent cross-collaboration between experts indicates a certain degree of unity and coherence in the field of mathematics. A newly discovered and surprising application of Theory A to Theory B is *not* a rare occurrence and is often put forward as the strongest evidence of the unity of mathematics.

On the other hand, professional specialism is rife and there is an often unconscious snobbism that downgrades interests other than one's own. "Life is too short", one says, "for me to care about everything let alone to know everything." A consequence is that one professional may have precious little to chat about with his office neighbour. The fact that Brown University has two Departments of Mathematics, one pure and one applied, with an historic record of intra-professional isolation is significant evidence of disunity.

**Diachronic and cross cultural disunity**
Written mathematics is easily 4000 years old. It has been created by people and has served people for a variety of purposes. A mathematician lives in a sub-culture at a certain time and place. A piece of mathematics does not exist only in a sequence of special symbols because the naked symbols are essentially uninterpretable. The symbols are embedded in a cloud of knowledge, meanings, associations, experiences and imaginations that derive from the particularities of time, place, person and the enveloping society.

Pythagoras asserted that 3 is the first male number. In certain Christian theologies it is the number of the Godhead. If in Chaldean numerology the numbers 1, 10, 19 and 28 are "ruled by the sun", the meaning of and the belief in those words may easily escape our readers. 666 is the "Number of the Beast" and there are numerous other mystic and arithmetic relationships. One might now say that such statements as "3 is the first male number" are not mathematics, pure or applied, but such a judgment delimits the notion of mathematics to only what is acceptable to today's establishment.

More generally, there are mathematical narratives that attempt to interpret very old texts. Methodologies of interpretation go by the fancy name *hermeneutics* and a frequently utilised interpretive device goes by the name *present-centred history*. This is where the past is described by an interpreter who uses the full knowledge of developments of subsequent importance. But can the interpreter do better than that and really get into the minds of the past creators? Only imperfectly. Eleanor Robson, noted for her studies of ancient Babylonian mathematics, suggests as much. Nonetheless, efforts have been made in this direction and Robson refers in [12] to the work of Henk Bos and Herbert Mehrtens who considered the relation between mathematics and the enveloping society and to the work of David Bloor for whom mathematics was a social construction toute courte.

**Semantic ambiguity**
We may write down the sequence of symbols x⌐∩σΣ≡6 and claim that this is a piece of mathematics. But this claim cannot be substantiated on the basis of the symbols alone. To provide meaning, every mathematical statement must be embedded in a sentence written in some natural language (English, German, etc.). Furthermore, significance of these symbols as mathematics cannot be established if its knowledge were limited as a private revelation to one and only one person.

Here is an example of private jottings. Paul Valéry (1871–1945) was a distinguished French poet and litterateur. He was also a bit of a "math buff".

Here is a clip from a page of his *Cahiers* (Notebooks).

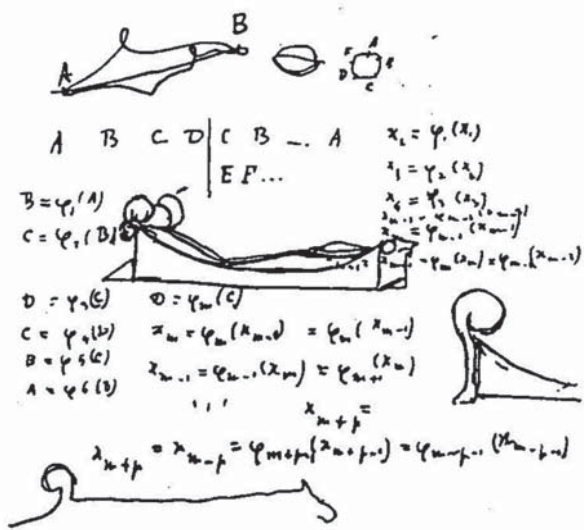

**Jottings extracted from Paul Valéry, *Cahiers II (1900–1902)*, Paris 1957, p. 608. Copyright by Centre National de la Recherche Scientifique.**

What is going on here? For a discussion of the attempt to interpret these jottings, see P.J.D.: *Bridging Two Cultures: Paul Valéry*.[2]

**Semiotic ambiguity and the problem of synonymy**
Can it be determined when two mathematical statements, phrased differently, are asserting the same thing? "The awkwardness of equality" has been discussed by Barry Mazur in his article *When is one thing equal to some other thing*?[3]

---

[2] *Svenska Matematikersamfundet Medlemsutskicket*, 15 May 2009. Reprinted in: *The Best Writing on Mathematics*. Mircea Pitici, 2010.
[3] Pp. 221–241 of [9].

**Scepticism, doubts, non-acceptance**

While each of the sections of this article could be elaborated into a book, the doubts, dilemmas and angst surrounding certain mathematical concepts could be shaped into a mathematical comedy. They have been treated in extenso by a variety of authors.[4] We shall recall briefly and simplistically and not in historical order a number of the objections that have been raised.

- Zero: How can zero be a number? Common sense asks: how can nothing be treated as something?
- One: A number must express numerosity or multiplicity.
- Fractions: How can 1/2 = 2/4 when half a pie is not the same as two pieces of pie that have been cut into quarters?
- Negative numbers: How can less than nothing be something? Scepticism about the interpretation of negatives lasted far into the 19th century. António José Teixeira, mathematics professor in Coimbra, wrote in 1890 that he did not like the proportion $1:-1::-1:1$ and asserted that "the negative quantities do not possess any arithmetical existence".
- Irrational numbers: $\sqrt{2}$ derives from the diagonal of the unit square. It is a line segment whose length exists and yet has no existence as a number. (Note the legend of Pythagoras offering up a hecatomb of oxen upon this discovery.)
- Imaginary numbers: $\sqrt{-1}$, etc. "Imaginary numbers are a fine and wonderful refuge of the divine spirit almost an amphibian between being and non-being" – Leibniz.
- Infinity: How can there be an infinity when the concept is germinated, represented or defined and operated on by a finite number of symbols? Infinity, so some claim, is a property only of the Godhead.
- Infinitesimals: These are "the ghosts of departed quantities" according to the philosopher George Berkeley.

The ideas surrounding number, particularly the idea of infinity, have been fertile fields for philosophers, theologians, neo-Platonists, numerologists, mystics and even cranks to plough. Indeed, in the early days of the development of the subject, mathematicians could hardly be distinguished from these other "occupations".

Moving briefly to mathematical functions (i.e. curves, graphs), we come across another objection:

- Dirac Function: How can a function that is zero on $(-\infty, +\infty)$ except at one point have a positive area "underneath" it?

Over the years – sometimes it took centuries – the concepts just mentioned and considered problematic have been totally absorbed into mainstream mathematics by, among other things, having been embedded within axiomatic, deductive formalisations of a traditional type. Moreover, and this is of great importance, these concepts have proved useful to science, technology and even to mathematics itself, as well as to a wide variety of humanistic concerns.

Thus, it would seem that both within and without mathematics, utility confers ontological reality and justificatory legitimacy. An increase in utility is accompanied by additional legitimacy and an abatement of scepticism. Yet the concept of utility and the ideas of "more" utility and "less" utility are admittedly vague and have been questioned. Utility to what and to whom? The late Richard Hamming, a very down-to-earth and practical type (he has been dubbed a techno-realist[5]), referring to one standard and basic concept of what is now termed "real analysis", wrote:

> "Does anyone believe that the difference between the Lebesgue and Riemann integrals can have physical significance, and that whether say, an airplane would or would not fly could depend on this difference? If such were claimed, I should not care to fly in that plane."

Profitable utility as a pre-condition for legitimacy stands in low regard in certain portions of the mathematical community. Wasn't it Euclid who is quoted as saying "Give the student a coin for he demands to profit from what he learns"? Other criteria besides utility exist for justification, legitimacy and acceptance: for example, the process of mathematical proof or physical verification. There are other instances of non-deductive approaches and all of these add up to what might be called modes of theorematic evidence.

Zermelo did not believe Gödel's proof. The scepticism of Kronecker, Poincaré, Zermelo, E. Picard, Brouwer, Hermann Weyl, Wittgenstein, Errett Bishopp, etc. regarding the concepts of Cantor has been well documented. Another quote from Richard Hamming sums this up:

> "I know that the great Hilbert said 'We shall not be driven out of the paradise that Cantor has created for us', and I reply I see no reason for walking in."

In the opposite direction, the great G. H. Hardy asserted he was interested in mathematics only as a creative art and that none of his work in it was of any utility.

**Philosophic ambiguity**

Originally there was one formal geometry: that of Euclid. After Bolyai and Lobatchevski, there were three and after Riemann, an infinity of geometries. Prior to the end of the 19th century, there was one philosophy of mathematics: that of Platonism. Now there are easily a number of distinguishable philosophies and together with these there are variations that exhibit the Freudian "narcissism of slight differences". There are platonism, logicism, formalism, intuitionism (constructivism), empiricism, conventionalism and culturalism/experientialism. In probability theory there are the frequentists vs. the interpretations of likelihood. We call ourselves "mild social constructivists". Others have called us simply "mavericks".

---

[4] E.g. Brian Rothman, *Signifying Nothing: The Semiotics of Zero*, or Gert Schubring, *Conflicts Between Generalization, Rigor, and Intuition: Number Concepts Underlying the Development of Analysis in 17–19th Century France and Germany*. See also: P.J. Davis: *Mathematics and Common Sense*, Chap. 12.

[5] *Richard Hamming: An Outspoken Techno-Realist*. SIAM NEWS, 16 Nov 1998.

Despite the often cited saying of Pythagoras that "All is Number", we doubt whether the Pythagoreans would have called the Songs of the Sirens mathematics. There is a separate and unique corpus of material, created over the centuries by the human intellect, that has been gathered together and has been called mathematics. Unity vs. disunitry? *E pluribus unum*.

**Acknowledgements**

We wish to acknowledge fruitful discussions with Sussi Booß-Bavnbek, Ernest S. Davis and Katalin Munkacsy.

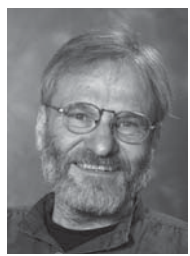

*Bernhelm Booß-Bavnbek [booss@ruc.dk] is a senior lecturer of global analysis and mathematical modelling at Roskilde University (Denmark). His newest book "Index Theory with Applications to Mathematics and Physics" is in print.*

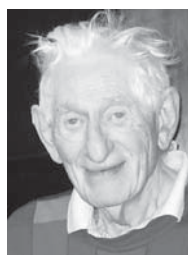

*Philip J. Davis [philip_davis@brown.edu] is a professor emeritus of applied mathematics at Brown University and is an independent writer, scholar and lecturer. He lives in Providence, Rhode Island, USA.*

# Episciences: a publishing platform for Open Archive Overlay Journals

Jean-Pierre Demailly (Université Joseph Fourier, Saint-Martin d'Hères, France) on behalf of the epi-math committee

While ever-increasing prices of scholarly publications have raised concerns for a long time, a growing number of interesting evolutions have taken place in the scientific community, especially in mathematics: a universally accepted document format (TeX/LaTeX), a widely adopted electronic open archive with worldwide coverage (arXiv) and more recent attempts to create open access journals and open discussion forums. In the context of commercial publishers seemingly unlikely to propose affordable and sufficiently open solutions for scientific publications, the mathematical community seems ready to adopt new publishing models and get actively involved in the related developments. Since creating open access electronic platforms has become technically easy, the obstacles are more historical and sociological than anything else. What is needed is a convenient framework that gives our community the tools to assess, correct, certify, archive and make widely available their production.

"Episciences" is a project hosted by the CCSD (Centre pour la Communication Scientifique Directe, located at Lyon University and attached to CNRS/IN2P3, a network of nuclear physics laboratories), in collaboration with Institut Fourier (a mathematics research department at Grenoble University, France). The CCSD develops the open archive "HAL" and also maintains a complete mirror of arXiv and an interface between HAL and arXiv. The aim of "Episciences" is to provide a publishing platform that makes it as easy as possible to host, run or create open archive overlay journals (hopefully a large number of them, independent of each other). Technically, the platform will rely on the HAL archive.

An overlay journal is a scientific journal that is focused on the peer-review process, and backed by one or several open archives for its diffusion and data handling. The first goal is to make it possible to efficiently run a journal at minimal cost and the second goal is to enforce a unified open access to the electronic version. In recent years, some print journals briefly became overlay journals but this experiment could not be sustained over an extended period. Episciences is an attempt at electronic-